%%%%%%%%%%%%%%%%%%%%%%%%%%%%
\input amstex.tex
\documentstyle{amsppt}
\magnification=1200
\loadbold
\def\Co#1{{\Cal O}_{#1}}
\def\Fi#1{\Phi_{|#1|}}
\def\fei#1{\phi_{#1}}

\def\lrw{\longrightarrow}

\def\Bbbp1{{\Bbb P}^1}
\def\simlin{\sim_{\text{\rm lin}}}
\def\simnum{\sim_{\text{\rm num}}}
\def\Div{\text{\rm Div}}

\def\dim{\text{\rm dim}}
\def\roundup#1{\ulcorner{#1}\urcorner}
\def\llrcorner#1{\llcorner{#1}\lrcorner}
\def\Bbbq{\Bbb Q}
\def\kod{\text{\rm kod}}
\TagsOnRight
%%%%%%%%%%%%%%%%%%%%%%%%%%%%%%%%%%%%%%%%%%%%%%%%%%
\topmatter
\title
Canonical stability in terms of singularity index for algebraic
threefolds
\endtitle
\author   {\smc By  Meng CHEN}
\endauthor
\affil
{\it Department of Applied Mathematics, Tongji University\\
 Shanghai, 200092, China\\
E-mail: {\bf chenmb\@online.sh.cn}} 
\endaffil
\thanks *
The author was partially supported by both the Post-Doc Fellowship of the Mathematisches Institut der Universit$\ddot{\text{a}}$t G$\ddot{\text{o}}$ttingen, Germany and the National Natural Science Foundation 
of China.
\endthanks
\endtopmatter
%%%%%%%%%%%%%%%%%%%%%%%%%%%%%%%%%%%%%%%%%%%%%%%%%
\document
%%%%%%%%%%%%%%%%%%%%%%%
%\baselineskip 16pt
\rightheadtext{Canonical stability for algebraic 3-folds}
\leftheadtext{M. Chen}
\centerline{
{\it (Received 4 February 2000; revised 19 April 2000)
\footnote{Accepted for publication in {\bf Math. Proc. Camb. Phil. Soc.}}
}}
%%%%%%%%%%

\head {\rm Introduction}\endhead
%%%%%%%%%%%

Throughout the ground field is always supposed to be 
algebraically closed of characteristic zero.
Let $X$ be a smooth projective threefold of general type, denote by
$\fei{m}$
the m-canonical map of $X$ which is nothing but the rational map naturally
associated
with the complete linear system $|mK_X|$. Since, once given such a 3-fold
$X$,
$\fei{m}$ is birational whenever $m\gg 0$, thus a quite interesting thing is
to find the optimal bound for such an $m$.
This bound is important because it is not only crucial to the classification
theory,
but also strongly related to other problems. For example, it can be applied
to
determine the order of the birational automorphism group of $X$ (\cite{21},
Remark in
\S1). To fix the terminology, we say
that $\fei{m}$ is {\it stably birational} if $\fei{t}$ is birational onto
its
image for all $t\ge m$.
It is well-known that the parallel problem in surface case was solved by
Bombieri (\cite{1}) and others. In the 3-dimensional case,  many authors have
studied the problem, in quite different ways. Because, in this paper, we are interested
in the results obtained by M. Hanamura (\cite{7}), we don't plan to
mention more references here. According to 3-dimensional MMP, $X$ has a
minimal model which is a normal projective 3-fold with only ${\Bbb
Q}$-factorial
terminal
singularities. Though $X$ may have many minimal models, the singularity
index
(namely the canonical index) of any of its minimal models is uniquely
determined by
$X$. Denote by $r$ the canonical index of minimal models of $X$.
When $r=1$, we know that
$\fei{6}$
is stably birational by virtue of \cite{3}, \cite{6}, \cite{13} and
\cite{14}.
When $r\ge 2$, M. Hanamura proved the following theorem.

\proclaim{\smc Theorem 0} (Theorem (3.4) of \cite{7})
Let $X$ be a smooth projective threefold of general type
with a minimal model of the canonical index $r$. Then $\fei{n_0(r)}$ is
stably birational onto its
image, where $n_0(r)$ is a function defined as 
$$\alignat4
 \hfill\ \ &\ \ \ \ * \ \ &\ \ r=2\ \ &\ \ 3\le r\le 5\ \ &\ \ r\ge 6\hskip3.5cm\\
\hfill\ \ &\ \ n_0(r)\ \ &\ \ 13\ \ \ \ &\ \ \ 4r+4\ \ &\ \ \ \ 4r+3.\hskip3.5cm
\endalignat$$
\endproclaim

Noting that the output $n_0(1)$ of Hanamura's method is actually $7$ (rather
than  6), it is reasonable to believe that the bound in Theorem 0 is not
optimal.
On the other hand, we don't know whether the canonical index $r$ is bounded
or
not, actually $r$ can be strangely large for some $X$.
This suggests that to find the optimal bounds for $n_0(r)$ should still be
very
interesting.
As far as our method can tell in this paper, the results are as the following

\proclaim{\smc Main Theorem} Let $X$ be a smooth projective threefold of general
type
with a minimal model of the canonical index $r$. Then

(i) $\fei{m}$ is generically finite whenever $m\ge l_0(r)$, where $l_0(r)$
is a function defined as   
$$\alignat4
 \hfill\ \ &\ \ \ \ * \ \ &\ \ r=2\ \ &\ \ 3\le r\le 5\ \ &\ \ r\ge 6\hskip3.5cm\\
\hfill\ \ &\ \ l_0(r)\ \ &\ \ 10\ \ \ \ &\ \ \ 2r+5\ \ &\ \ \ \ 2r+4.\hskip3.5cm
\endalignat$$

(ii) $\fei{m_0(r)}$ is stably birational onto its image, where $m_0(r)$ is a
function defined as  
$$\alignat7
 \hfill\ \ &\ \ \ \ * \ \ &\ \ r=2\ \ &\ \ r=3\ \ &\ \ r=4\ \ &\ \ r=5\ \ &\ \ r=6\ \ &\ \ r\ge 7\hskip2cm\\
\hfill\ \ &\ \ m_0(r)\ \ &\ \ 11\ \ \ \ &\ \ \ \ 15\ \ &\ \ 17\ \ \ \ &\ \ \ \ 18\ \ &\ \ 19\ \ \ \ &\ \ 2r+6.\hskip2cm
\endalignat$$
\endproclaim

As an application of our method, we shall present the following

\proclaim{\smc Corollary} Let $X$ be a smooth projective 3-fold of general
type.
Then $\fei{9}$ is birational if $p_g(X)\ge 2$. 
\endproclaim

\remark{Remark} The above 
corollary is an improvement to Koll\'ar's result (Corollary 4.8 of
\cite{11}) that $\fei{16}$ is birational if $p_g(X)\ge 2$.
Actually, Koll\'ar proved there that $\fei{11k+5}$ is birational if
$P_k:=h^0(X, kK_X)\ge 2,$
 where $k$ is a positive integer. Recently, \cite{4}
improved this result to the level that either $\fei{7k+3}$ or $\fei{7k+5}$
is birational under the same condition.
\endremark

For readers' convenience, we briefly explain the whole technique of this
paper.
According to Hanamura's result that $|(r+2)K_X|$ is not composed of a
pencil,
we can take a general member $S_2$ of the movable part of this system.
Actually
we can suppose that $S_2$ is smooth.
Then we use the Matsuki-Tankeev principle to reduce the birationality
problem to a parallel one for the adjoint system $|K+L|$ on the surface
$S_2$ which is a smooth projective surface of general type. We shall
inevitably
treat a very delicate case in which  $L$ is the round-up of certain nef and big 
${\Bbb Q}$-divisor $A$, i.e. $L=\roundup{A}$. Instead of applying Reider's result, we go on reducing to the problem on a curve. The technical point is to estimate the degree of the divisor in question on the curve. The Kawamata-Ramanujam-Viehweg vanishing theorem played an important role in the whole context. 

\head {\rm 1. Preliminaries} \endhead

Let $X$ be a normal projective variety of dimension $d$. We denote by
$\Div(X)$
the group of Weil divisors on $X$. An element $D\in\Div(X)\otimes{\Bbb Q}$
is called a ${\Bbb Q}$-{\it divisor}. A ${\Bbb Q}$-divisor $D$ is
said to be ${\Bbb Q}$-{\it Cartier} if $mD$ is a Cartier divisor for some
positive integer $m$. For a ${\Bbb Q}$-Cartier divisor
$D$ and an irreducible curve $C\subset X$, we can define the intersection
number $D\cdot C$ in a natural way. A $\Bbbq$-Cartier divisor $D$ is called
{\it nef} (or {\it numerically effective}) if $D\cdot C\ge 0$ for any
effective curve
$C \subset X$. A nef divisor $D$ is called {\it big} if $D^d>0$. We say that
$X$ is $\Bbbq$-{\it factorial} if every Weil divisor on $X$ is
$\Bbbq$-Cartier.
For a Weil divisor $D$ on $X$, write $\Co{X}(D)$ as the corresponding
reflexive
sheaf.
Denote by $K_X$ a canonical divisor of $X$, which is a Weil divisor. $X$ is
called
{\it minimal} if $K_X$ is a nef $\Bbbq$-Cartier divisor. $X$ is said to be
of general
type if the Kodaira dimension $\kod(X)=\dim(X)$. For a positive integer $m$, we set
$\omega_X^{[m]}:=\Co{X}(mK_X)$
and call
$P_m(X):=\dim_{k}H^0(X, \omega_X^{[m]})$
the {\it m-th plurigenus} of
$X$. We remark that $P_m(X)$ is an important birational invariant.

$X$ is said to have only {\it canonical singularities} (resp. {\it terminal
singularities}) according to Reid (\cite{15}) if the following two conditions
hold:

(i) for some positive integer $r$, $rK_X$ is Cartier;

(ii) for some resolution $f:Y\lrw X$,
$K_Y=f^*(K_X)+\sum a_iE_i$
for $0\le a_i\in \Bbbq$ (resp. $0<a_i$) $\forall i$, where the $E_i$ vary
amongst all the
exceptional divisors on $Y$.
The minimal $r$ that satisfies (i) is called {\it the canonical index}
of $X$ and can be also written as $r(X)$.

According to Mori's MMP (\cite{10}, \cite{12}), when $V$ is a smooth
projective threefold,
there exists a birational map $\sigma: V\dashrightarrow X$ where $X$ can be
a
minimal 3-fold with only $\Bbbq$-factorial terminal singularities. Usually,
$X$ is not uniquely determined by $V$, but the canonical index $r(X)$ is.

Let $D=\sum a_iD_i$ be a $\Bbbq$-divisor on $X$ where the $D_i$ are distinct
prime divisors and $a_i\in\Bbbq$. We define
$$\align
&\text{the round-down}\ \llrcorner{D}:=\sum\llrcorner{a_i}D_i,\
\text{where}\
\llrcorner{a_i}\ \text{is the integral part of}\ a_i,\\
&\text{the round-up}\  \roundup{D}:=-\llrcorner{-D},\\
&\text{the fractional part}\ \{D\}:=D-\llrcorner{D}.
\endalign$$

\remark{Remark 1.1} Suppose that $X$ has only canonical singularities and
that $f:V\lrw X$ is a resolution. We have
$$P_m(X)=h^0\Bigl(V,\Co{V}\bigl(\llrcorner{f^*(mK_X)}\bigr)\Bigr)
=h^0\Bigl(V,\Co{V}\bigl(\roundup{f^*(mK_X)}\bigr)\Bigr)
=P_m(V)$$
for any positive integer $m$.
\endremark

We always use the Kawamata-Ramanujam-Viehweg vanishing theorem in the
following
form.
\proclaim{\smc Vanishing Theorem} (\cite{9} or \cite{18}) Let $X$ be a smooth
complete
variety, $D\in\Div(X)\otimes\Bbbq$. Assume the following two conditions:

(i) $D$ is nef and big;

(ii) the fractional part of $D$ has supports with only normal crossings.
\noindent Then $H^i(X, \Co{X}(K_X+\roundup{D}))=0$ for all $i>0$.
\endproclaim

Most of our notations are standard within algebraic geometry except the
following
which we are in favor of:
$\simlin$ means {\it linear equivalence} while $\simnum$ means {\it
numerical equivalence}.

\head {\rm 2. Some lemmas} \endhead

\subhead 2.1 The Matsuki-Tankeev principle\endsubhead
This principle is tacitly used throughout our argument. Suppose $X$ is a
smooth
variety, $|M|$ is a base point free system on $X$ and $D$ is a divisor with
$|D|\ne\varnothing$. We want to know when $\Fi{D+M}$ is birational. The
following
principles are due to Tankeev and Matsuki respectively.

\proclaim{(P1)} (Lemma 2 of \cite{17}) Suppose $|M|$ is not composed of a
pencil, i.e.
$$\dim\Fi{M}(X)\ge 2$$ 
and take a general member $Y\in |M|$. If the restriction
of $\Fi{D+M}$ to $Y$ is birational, then $\Fi{D+M}$ is birational.
\endproclaim

\proclaim{(P2)} (see the proof of the main theorem in \cite{14}) Suppose $|M|$ is composed of a pencil and
take the Stein factorization of
$$\Fi{M}: X\overset{f}\to\lrw C\lrw W\subset {\Bbb P}^N,$$
where $W$ is the image of $X$ through $\Fi{M}$ and $f$ is a fibration onto a
smooth curve $C$. Let $F$ be a general fiber of $f$. If we know (say
by the vanishing
theorem) that $\Fi{D+M}$ can distinguish general fibers of $f$ (i.e. separates 
any two points of the respective fibers)
and
its restriction to $F$ is birational, then $\Fi{D+M}$ is also birational.
\endproclaim
%%%%%%%%%%%%%
\proclaim{\smc Lemma 2.2} Let $X$ be a smooth projective variety of dimension
$d$,
$D\in \Div(X)\otimes\Bbbq$ be a $\Bbbq$-divisor on $X$. Then the following
assertions are true:

(i) if $S$ is a smooth reduced irreducible divisor on $X$ and $S$ is not a fractional component of $D$, then
$\roundup{D}|_S\ge \roundup{D|_S}$;

(ii) if $\pi:X'\lrw X$ is a birational morphism, then
$\pi^*(\roundup{D})\ge
\roundup{\pi^*(D)}$.
\endproclaim
\demo{Proof}
This lemma is very easy to check.
\qed
\enddemo
\proclaim{\smc Lemma 2.3} Let $S$ be a smooth projective surface of general type
and
$L$ be a nef and big divisor on $S$. Then $\Fi{K_S+mL}$ is birational in 
the following cases:

(i) $m\ge 4$;

(ii) $m=3$ and $L^2\ge 2$.
\endproclaim
\demo{Proof} This is a direct result of Corollary 2 of \cite{16}.
\qed
\enddemo
\proclaim{\smc Lemma 2.4} (Lemma (3.2) of \cite{7}) Let $X$ be a minimal
threefold
of general type with canonical index $r\ge 2$. Then
$\dim\fei{mr+s}(X)\ge 2$
in the following cases, where $m$ is a positive integer and $0\le
s<r$:

(i) $r=2$ and $m\ge 3$;

(ii) $r=3$ and $m\ge 2$;

(iii) $r=4,\ 5$, $0\le s\le 2$ and $m\ge 2$; $r=4,\ 5$, $s\ge 3$ and $m\ge
1$;

(iv) $r\ge 6$, $0\le s\le 1$ and $m\ge 2$; $r\ge 6$, $s\ge 2$ and $m\ge 1$.
\endproclaim
\proclaim{\smc Lemma 2.5} Under the same assumption as in Lemma 2.4, the
plurigenus \newline
$P_{mr+s}(X)\ge 3$ 
in one of the following cases:

(i) $r=2$ and $m\ge 2$;

(ii) $r\ge 3$, $0\le s\le 1$ and $m\ge 2$; $r\ge 3$, $s\ge 2$ and $m\ge 1$.
\endproclaim
\demo{Proof}
This is obvious from the proof of Lemma (3.2) in \cite{7}. 
In order to be precise and to cite it 
many times in this paper, let us recall the estimation
there.
%%%%%%%%%%%

When $r\ge 3$, $r$ is even and $s\ge 2$, we have
$$\align
P_{mr+s}(X)&\ge
\frac{1}{12}\biggl\{2r^2m^3+(6s-3)rm^2+\Bigl(6s^2-6s-\frac{1}{2}r^2\Bigr)m
\biggr\}
(rK_X^3)\\
&\ge \frac{1}{8}(r^2+6r+9) \tag{2.1}
\endalign$$
%%%%%%%%%%%

When $r\ge 3$, $r$ is odd and $s\ge 2$, we have
$$\align
P_{mr+s}(X)&\ge
\frac{1}{12}\biggl\{(mr+s)(mr+s-1)(2mr+2s-1)+m\Bigl(-\frac{1}{2}r^3+
 \frac{1}{2}r\Bigr)  \\
&\ \ \ -s(s-1)(2s-1)\biggr\}(K_X^3) \\
 &\ge \frac{1}{8}(r^2+6r+8) \tag{2.2}
\endalign$$
%%%%%%%%%%%%

When $s=1$, we have
$$P_{mr+1}(X)\ge \frac{1}{12}r(m^2-1)(2rm+3)(rK_X^3) \tag{2.3}$$
%%%%%%%%%%%%

When $s=0$, we have
$$P_{mr}(X)\ge \frac{1}{12}r(m^2-1)(2rm-3)(rK_X^3) \tag{2.4}$$
The above four formulae give the result.
\qed
\enddemo
%%%%%%%%%%

The following lemma is sufficient to derive our result, though it seems that one
might exploit its potential.
%%%%%%%%%
\proclaim{\smc Lemma 2.6} Let $S$ be a smooth projective surface  of general type
and $L$ be a nef
divisor on $S$ such that $|L|$ gives a generically finite map. Then

(i) $L^2\ge h^0(S,L)-2$; if $\Fi{L}$ is not birational, then
$L^2\ge 2h^0(S, L)-4$.

(ii) if $p_g(S)>0$, then $L^2\ge 2h^0(S, L)-4$.

(iii) $K_S+L$ is always effective.
\endproclaim
\demo{Proof}
The first part is trivial. One should note that a non-degenerate surface in
${\Bbb P}^N$
has degree $\ge N-1$. In order to prove the second part, we may suppose
that $|L|$ is base point free. Let $C$ be a general member of $|L|$, then
$$h^0(C, L|_C)\ge h^0(S, L)-1.$$
Noting that $C$ is moving and that we have the following exact sequence
$$0\lrw\Co{S}(K_S-C)\lrw\Co{S}(K_S)\lrw \Co{C}(K_S|_C)\lrw 0,$$
the inclusion
$$H^0(S, K_S-C)\hookrightarrow H^0(S, K_S)$$
is proper. So
$h^0(C, K_S|_C)>0$, which means $h^1(C,L|_C)>0$ and $L|_C$ is special. Thus
 Clifford's theorem implies that
$$L^2=deg(L|_C)\ge 2h^0(C, L|_C)-2\ge 2h^0(S, L)-4.$$
Finally, (iii) is an easy exercise by Riemann-Roch.
\qed
\enddemo
%%%%%%%%%%%%

\proclaim{\smc Lemma 2.7}
Let $X$ be a smooth projective variety of dimension $\ge 2$. Let $D$ be a divisor on $X$, $h^0(X, \Co{X}(D))\ge 2$ and $S$ be a smooth irreducible divisor on $X$ such that $S$ is not a fixed component of $|D|$. Denote by $M$ the movable part of $|D|$ and by $N$ the movable part of $|D|_S|$ on $S$. Suppose the natural restriction map
$$H^0(X, \Co{X}(D))\overset\theta\to\lrw H^0(S, \Co{S}(D|_S))$$
is surjective. Then $M|_S\ge N$ and thus
$$h^0(S,\Co{S}(M|_S))=h^0(S,\Co{S}(N))=h^0(S,\Co{S}(D|_S)).$$
\endproclaim
\demo{Proof} Denote by $Z$ the fixed part of $|D|$. Because $S$ is not a fixed component of $|Z|$, we see that $Z|_S\ge 0$. Thus $D|_S\ge M|_S$. Considering 
the natural map
$$H^0(X,\Co{X}(M))\overset\theta_0\to\lrw H^0(S,\Co{S}(M|_S)),$$
we have
$$h^0(S,\Co{S}(M|_S))\ge\dim_{\Bbb C}im(\theta_0)=\dim_{\Bbb C}im(\theta)
=h^0(S,\Co{S}(D|_S)).$$
This means $h^0(S,\Co{S}(M|_S))=h^0(S,\Co{S}(D|_S))$ and so $M|_S\ge N$ on $S$.
\qed\enddemo

%%%%%%%%%%%%
\head {\rm 3. The generic finiteness} \endhead
%%%%%%%%%
This section is devoted to study the generic finiteness of $\fei{m}$.
Whenever
we mention a minimal 3-fold, we mean one with only $\Bbbq$-factorial terminal
singularities. The following theorem is the easy part of the Main Theorem.
%%%%%%%%%%%
\proclaim{\smc Theorem 3.1} Let $X$ be a projective minimal 3-fold of general
type
with $2\le r(X)\le 5$. Then $\fei{4r+3}$ is stably birational.
\endproclaim
\demo{Proof} According to Lemma 2.5, $P_{m_1}(X)\ge 3$
for $m_1\ge r+2$. Take necessary blowing-ups $\pi:X'\lrw X$
along nonsingular centers, according to Hironaka, such that $X'$ is smooth
and
$|m_1K_{X'}|$ defines a morphism (of course, $|m_1K_{X'}|$ may have fixed
components).
Denote by $M$ the movable part of $|m_1K_{X'}|$. We have
$$|K_{X'}+3\pi^*(rK_X)+M|\subset |(m_1+3r+1)K_{X'}|.$$
First we note that $K_{X'}+3\pi^*(rK_X)$ is effective according to Lemma
2.5.
If $|M|$ is not composed of a pencil, then a general member of it is
an irreducible smooth projective
surface $S$ of general type. Set $L:=\pi^*(rK_X)|_S$, which is a nef and big
Cartier divisor on $S$ with $L^2\ge 2$. Using the vanishing theorem, we get
$$|K_{X'}+3\pi^*(rK_X)+S|\bigm|_S=|K_S+3L|.$$
The right system gives a birational map by Lemma 2.3. So (P1) implies what
we want
in this case.
If $|M|$ is composed of a pencil, we take the Stein-factorization of
$$\Fi{M}:X'\overset{f}\to\lrw C\lrw W,$$
where $W$ is the image of $X'$ through
$\Fi{M}$ and $f$ is a fibration onto the smooth curve $C$. Generically, $M$
can be
written as a disjoint union of fibers of $f$, i.e.\newline
$M\simlin\sum_{i=1}^aF_i.$
The $F_i$'s are irreducible smooth surfaces of general type. The effectiveness
of
$K_{X'}+3\pi^*(rK_X)$ implies that
$|K_{X'}+3\pi^*(rK_X)+M|$ 
can distinguish
general fibers of the morphism $\Fi{M}$. On the other hand, we have the
following
exact sequence
$$H^0\Bigl(X', K_{X'}+3\pi^*(rK_X)+M\Bigr)\lrw\oplus_{i=1}^a H^0(F_i,
K_{F_i}+3L_i)\lrw 0$$
where $L_i:=\pi^*(rK_X)|_{F_i}$ is a nef and big divisor with $L_i^2\ge 2$.
This shows that the system
$|K_{X'}+3\pi^*(rK_X)+M|$
can also distinguish different
components in a general fiber of $\Fi{M}$ and the restriction to each $F_i$
gives
a birational map. Thus, by (P2), we have completed the proof.
\qed
\enddemo
%%%%%%%%%%%%%%%
\remark{Remark 3.2}
Throughout this paper, we shall deal with the same situation as in the proof
of Theorem 3.1. In order to avoid unnecessary redundancy, we give the
definition of so-called {\it generic irreducible element} of a moving system
$|M|$ on a
variety $V$. Using our notations in the proof of Theorem 3.1, we shall
call $S$
(respectively, $F_i$) a generic irreducible element of $|M|$ ignoring
whether
it is composed of a pencil or not. In our case, we always use both (P1) and
(P2).
\endremark
%%%%%%%%%%%%%%%

\proclaim{\smc Theorem 3.3}
Let $X$ be a projective minimal threefold of general type with the 
canonical index $r$. Then $\fei{m}$ is generically finite whenever
$m\ge l_0(r)$, where $l_0(r)$ is a function defined as the following
$$l_0(r)=\cases 10\ \ &\text{if}\ r=2,\\
2r+5\ \ &\text{if}\ 3\le r\le 5,\\
2r+4\ \ &\text{if}\ r\ge 6.\endcases$$
\endproclaim
\demo{Proof} The idea of the proof is quite simple. We formulate our proof through steps.

Step 1. Set up for the proof.

First, we define
$$m_2=\cases 6\ \ &\text{if}\ r=2,\\
r+3\ \ &\text{if}\ 3\le r\le 5,\\
r+2\ \ &\text{if}\ r\ge 6.\endcases$$
Take a birational modification $\pi:X'\lrw X$, according to Hironaka, such
that

(1) $X'$ is smooth;

(2) $|m_2K_{X'}|$ defines a morphism;

(3) the fractional part of $\pi^*(K_X)$ has supports with only normal
crossings.
%%%%%%%%%%%%%

Denote by $M_2'$ the movable part of $|m_2K_{X'}|$ and by $S_2'$ a general
member
of $|M_2'|$. {}From Lemma 2.4, we know that $\dim\fei{m_2}(X)\ge 2$.
Given any integer $t\ge r+2$, we know that $|tK_{X'}|$ is always effective according to
Lemma 2.5. If $|M_2'|$ has already given a generically finite map, then
$\fei{t+m_2}$ is generically finite and thus the theorem is true in this
situation. So from now on, we suppose
$\dim\fei{m_2}(X)=2$. {}From (2.1), (2.2), (2.3) and (2.4), we have
$$h^0(X', S_2')=P_{m_2}(X)\ge
\cases 12\ \ &\text{if}\ r=2,\\
7\ \ &\text{if}\ r=3,\\
\frac{1}{8}(r^2+10r+24)\ \ &\text{if}\ 4\le r\le 5,\\
\frac{1}{8}(r^2+6r+8)\ \ &\text{if}\ r\ge 6.\endcases$$
On the surface $S_2'$, we set $L_2':=M_2'|_{S_2'}$. Then $|L_2'|$ is composed of
a pencil.
We can write
$$L_2'\simlin\sum_{i=1}^{a_2'}C_i\simnum a_2'C,$$
where
$a_2'\ge h^0(S_2', L_2')-1$
and $C$ is a generic irreducible element of $|L_2'|$.
Since
$h^0(S_2', L_2')\ge h^0(X', S_2')-1,$
we get
$$a_2'\ge
\cases 10\ \ &\text{if}\ r=2,\\
5\ \ &\text{if}\ r=3,\\
\frac{1}{8}(r^2+10r+8)\ \ &\text{if}\ 4\le r\le 5,\\
\frac{1}{8}(r^2+6r-8)\ \ &\text{if}\ r\ge 6.\endcases$$

Step 2. Reduce to the problem on a curve.

Because
$$|K_{X'}+\roundup{(t-1)\pi^*(K_X)}+S_2'|\subset |(t+m_2)K_{X'}|,$$
it is sufficient to prove the generic finiteness of
$\Fi{K_{X'}+\roundup{(t-1)\pi^*(K_X)}+S_2'}$. Noting that
$K_{X'}+\roundup{(t-1)\pi^*(K_X)}$
is effective, we only have to verify that
$$|K_{X'}+\roundup{(t-1)\pi^*(K_X)}+S_2'|\bigm|_{S_2'}$$
gives a generically finite map by virtue of (P1). 

The vanishing theorem gives
$$|K_{X'}+\roundup{(t-1)\pi^*(K_X)}+S_2'|\bigm|_{S_2'}
=|K_{S_2'}+\roundup{(t-1)\pi^*(K_X)}|_{S_2'}|,$$
so we are reduced to verify the same property for
$$|K_{S_2'}+\roundup{(t-1)\pi^*(K_X)}|_{S_2'}|.$$
Since
$$K_{S_2'}+\roundup{(t-1)\pi^*(K_X)}|_{S_2'}=
\Bigl(K_{X'}+\roundup{(t-1)\pi^*(K_X)}\Bigr)\bigm|_{S_2'}+L_2'$$
and $K_{X'}+\roundup{(t-1)\pi^*(K_X)}$ is effective by Lemma 2.5, 
the system 
$$|K_{S_2'}+\roundup{(t-1)\pi^*(K_X)}|_{S_2'}|$$
 can distinguish
general fibers of $\Fi{L_2'}$. So it is sufficient to show that
$$\Fi{K_{S_2'}+\roundup{(t-1)\pi^*(K_X)}|_{S_2'}}\bigm|_C$$
is a finite map for a generic irreducible element $C$ of $|L_2'|$.
%%%%%%%%%%%%%%%

Step 3. Verifying the finiteness on $C$.
%%%%%%%%%%%%

Since $m_2\pi^*(K_X)\ge S_2'$, we can write
$$m_2\pi^*(K_X)|_{S_2'}=L_2'+E_{\Bbb Q}
\simnum a_2'C+E_{\Bbb Q},$$
where $E_{\Bbb Q}$ is an effective $\Bbbq$-divisor on $S_2'$.
So we have
$$(t-1)\pi^*(K_X)|_{S_2'}\simnum\frac{(t-1)a_2'}{m_2}C+\frac{t-1}{m_2}E_{\Bbbq
}$$
and
$$(t-1)\pi^*(K_X)|_{S_2'}-C-\frac{1}{a_2'}E_{\Bbbq}\simnum
(t-1)\biggl(1-\frac{m_2}{(t-1)a_2'}\biggr)\pi^*(K_X)|_{S_2'}.$$
Set
$\alpha:=1-\frac{m_2}{(t-1)a_2'},$
it is easy to verify that $\alpha>0$. This shows that
$$H^1\Bigl(S_2',
K_{S_2'}+\roundup{(t-1)\pi^*(K_X)|_{S_2'}-\frac{1}{a_2'}E_{\Bbbq}}-C\Bigr)
=0$$
according to the vanishing theorem. Thus we have the exact sequence
$$H^0\Bigl(S_2',
K_{S_2'}+\roundup{(t-1)\pi^*(K_X)|_{S_2'}-\frac{1}{a_2'}E_{\Bbbq}}\Bigr)
\lrw H^0(C, K_C+D)\lrw 0,$$
where
$D:=\roundup{(t-1)\pi^*(K_X)|_{S_2'}-\frac{1}{a_2'}E_{\Bbbq}}\bigm|_C$
is a divisor on $C$ with positive degree. In fact,
$$\deg(D)\ge (t-1)\alpha\pi^*(K_X)|_{S_2'}\cdot C>0.$$
Noting that $C$ is a smooth curve of genus $\ge 2$, we see that $|K_C+D|$ gives a finite map.
Therefore
$$\Fi{K_{S_2'}+\roundup{(t-1)\pi^*(K_X)|_{S_2'}-\frac{1}{a_2'}E_{\Bbbq}}}
\bigm|_C$$
is generically finite. So
$\Fi{K_{S_2'}+\roundup{(t-1)\pi^*(K_X)}|_{S_2'}}\bigm|_C $
is also generically finite. This derives the generic finiteness of
$\fei{t+m_2}$.
\qed
\enddemo
%%%%%%%%%%%
\head  {\rm 4. The birationality:\  $r\ge 4$} \endhead
%%%%%%%%%%%

Suppose $X$ is a projective minimal 3-fold of general type with the
canonical
index $r$. First we take a birational modification  $\pi: X'\lrw X$,
according to Hironaka, such that

(i) $X'$ is smooth;

(ii) the system $|(r+2)K_{X'}|$ defines a morphism;

(iii) the fractional part of $\pi^*(K_X)$ has supports with only normal
crossings.

Denote by $M_2$ the movable part of $|(r+2)K_{X'}|$. We know from Lemma 2.4
that $|M_2|$ is not composed of a pencil when $r\ge 6$.

{}From now on, 
we assume that $r\ge 4$ and that $|(r+2)K_{X'}|$ is not composed of a
pencil.
Take a general member $S_2\in |M_2|$. Then $S_2$ is a smooth projective
surface of general type.
Set $L_2:=M_2|_{S_2}$. Then $L_2$ is a nef divisor on the surface
$S_2$.
We have already known from the proof of Lemma 2.5 that, for $r\ge 4$, 
$$h^0(X', S_2)=P_{r+2}(X)\ge \frac{1}{8}(r^2+6r+8).$$
%%%%%%%%%%

\proclaim{\smc Theorem 4.1} Let $X$ be a projective minimal 3-fold of general
type
with the canonical index $r\ge 4$. If $\dim\fei{r+2}(X)=3$, then
$\fei{m_1(r)}$ is stably birational where
$$m_1(r)=\cases 16, \ \ &\ \text{if} \ r=4\\
2r+7,\ \ &\ \text{if} \ 5\le r\le 6\\
2r+6,\ \ &\ \text{if} \ r\ge 7.
\endcases$$ 
\endproclaim
\demo{Proof} Given an integer $t_1>0$,
first we note that
$$|K_{X'}+\roundup{(t_1+r+2)\pi^*(K_X)}+S_2|\subset |(t_1+2r+5)K_{X'}|.$$
In order to prove the birationality of $\fei{t_1+2r+5}$, it is sufficient to
prove the same thing for
$\Fi{K_{X'}+\roundup{(t_1+r+2)\pi^*(K_X)}+S_2}.$

Step 1. Reduce to the problem on a surface.

Since
$K_{X'}+\roundup{(t_1+r+2)\pi^*(K_X)}$ is effective by Lemma 2.5, it is enough to verify the same
thing for its restriction to $S_2$ by virtue of (P1). The vanishing theorem
gives the exact sequence
$$\align
&H^0\bigl(X',K_{X'}+\roundup{(t_1+r+2)\pi^*(K_X)}+S_2\bigr)\lrw\\
&H^0\bigl(S_2, K_{S_2}+\roundup{(t_1+r+2)\pi^*(K_X)}|_{S_2}\bigr)\lrw 0.
\endalign$$
This means
$$|K_{X'}+\roundup{(t_1+r+2)\pi^*(K_X)}+S_2|\bigm|_{S_2}=
\bigm|K_{S_2}+\roundup{(t_1+r+2)\pi^*(K_X)}|_{S_2} \bigm|.$$
And from Lemma 2.2, we have
$$\bigm|K_{S_2}+\roundup{(t_1+r+2)\pi^*(K_X)|_{S_2}}\bigm|\subset
\bigm|K_{S_2}+\roundup{(t_1+r+2)\pi^*(K_X)}|_{S_2}\bigm|.$$
So, sometimes,  it is enough to show that
$\bigm|K_{S_2}+\roundup{(t_1+r+2)\pi^*(K_X)|_{S_2}}\bigm|$
gives a birational map. Under the assumption of this theorem, it is obvious
that $|L_2|$ gives a generically finite map.

Step 2. Reduce to the problem on a curve.

We suppose $\overline{C}$ is the general member of the movable part of $|L_2|$.
Since
$$K_{S_2}+\roundup{(t_1+r+2)\pi^*(K_X)}|_{S_2}=\bigl(K_{X'}+
\roundup{(t_1+r+2)\pi^*(K_X)}\bigr)|_{S_2}+L_2\ge L_2,$$
by (P1), we only have to verify the birationality of
$\Fi{K_{S_2}+\roundup{(t_1+r+2)\pi^*(K_X)}|_{S_2}}\bigm|_{\overline{C}}$
for a general member $\overline{C}$.
It is obvious that
$$\align
K_{S_2}+\roundup{(t_1+r+2)\pi^*(K_X)}|_{S_2}&\ge
K_{S_2}+\roundup{t_1\pi^*(K_X)}|_{S_2}+L_2\\
&\ge K_{S_2}+\roundup{t_1\pi^*(K_X)|_{S_2}}+\overline{C}.
\endalign $$
It's sufficient to show that $\Fi{K_{S_2}+\roundup{t_1\pi^*(K_X)|_{S_2}}+\overline{C}}|_{\overline{C}}$
is birational.

Step 3. Verifying the embedding on $\overline{C}$.

The vanishing theorem gives
$$\bigm|K_{S_2}+\roundup{t_1\pi^*(K_X)|_{S_2}}+\overline{C} \bigm|
\Bigm|_{\overline{C}}=
|K_{\overline{C}}+D_0|,$$
where
$D_0:=\roundup{t_1\pi^*(K_X)|_{S_2}}|_{\overline{C}}$
is a divisor on $\overline{C}$ with 
$$\deg(D_0)\ge t_1\pi^*(K_X)|_{S_2}\cdot \overline{C}.$$
It's clear that the theorem follows whenever $\deg(D_0)\ge 3$.
Although, in general, $\pi^*(K_X)|_{S_2}\cdot \overline{C}$ is a positive rational number, we can estimate it in this situation.

Note that, if $|\overline{C}|$ has already given a birational map, then 
so does $|K_{S_2}+\roundup{(t_1+r+2)\pi^*(K_X)}|_{S_2}|$ because 
$$K_{S_2}+\roundup{(t_1+r+2)\pi^*(K_X)}|_{S_2}\ge L_2.$$
So we may suppose that $|\overline{C}|$ gives a generically finite, non-birational map on the surface $S_2$. According to Lemma 2.6(i), we get
$$\overline{C}^2\ge 2h^0(S_2, \overline{C})- 4\ge \frac{1}{4}(r^2+6r)-4.$$
Thus we get
$$(r+2)\pi^*(K_X)|_{S_2}\cdot \overline{C}\ge \overline{C}^2\ge
\frac{1}{4}(r^2+6r)-4.$$
One can easily obtain
$$\pi^*(K_X)|_{S_2}\cdot \overline{C}\cases
\ge 1,\ &\ \text{if}\ r=4\\
\ge \frac{10}{7},\ &\ \text{if}\ r=5\\
\ge \frac{7}{4},\ &\ \text{if}\ r=6\\
 >2, \ &\ \text{if}\ r\ge 7.\endcases$$
Thus we can see that, whenever $t_1\ge m_1(r)-2r-5$, 
$$\deg(D_0)\ge \roundup{t_1\pi^*(K_X)|_{S_2}\cdot \overline{C}}\ge 3.$$
We have proved that $\fei{m_1(r)}$ is stably birational under assumption of the theorem.
\qed
\enddemo
%%%%%%%%%%%%%%
\proclaim{\smc Theorem 4.2}
 Let $X$ be a projective minimal 3-fold of general type
with the canonical index $r\ge 4$. If $\dim\fei{r+2}(X)=2$, then
$\fei{2r+6}$ is stably birational.
\endproclaim
\demo{Proof}
In this case, we note that $|L_2|$ is a base point free pencil on the
surface
$S_2$. We can write
$$ L_2\simlin \sum_{i=1}^{a_2}C_i\simnum a_2C,$$
where
$$a_2\ge h^0(S_2, L_2)-1\ge \frac{1}{8}(r^2+6r-8)$$
and $C$ denotes a generic
irreducible element of $|L_2|$. Given an integer $t_2>0$,
we want to prove the birationality of $\fei{t_2+2r+5}$. For the same reason as in the proof of Theorem 4.1, 
it is sufficient to verify the birationality of 
$\Fi{K_{S_2}+\roundup{(t_2+r+2)\pi^*(K_X)}|_{S_2}}.$

Step 1. Reduce to the problem on a curve.

Since $K_{X'}+\roundup{(t_2+r+2)\pi^*(K_X)}$ is effective by Lemma 2.5, 
we can see that
$$K_{S_2}+\roundup{(t_2+r+2)\pi^*(K_X)}|_{S_2}\ge L_2.$$
This shows that
$\Fi{K_{S_2}+\roundup{(t_2+r+2)\pi^*(K_X)}|_{S_2}}$
can distinguish different fibers of $\Fi{L_2}$. On the other hand, we have
$$(t_2+r+2)\pi^*(K_X)\ge t_2\pi^*(K_X)+S_2.$$
So we get
$$K_{S_2}+\roundup{(t_2+r+2)\pi^*(K_X)}|_{S_2}\ge
K_{S_2}+\roundup{t_2\pi^*(K_X)}|_{S_2}+L_2.$$
The vanishing theorem gives the following exact sequence
$$H^0\Bigl(S_2, K_{S_2}+\roundup{t_2\pi^*(K_X)|_{S_2}}+L_2\Bigr)
\lrw\oplus_{i=1}^{a_2}H^0(C_i, K_{C_i}+D_i)\lrw 0$$
where the $D_i's$ are divisors on the curve $C_i$ with positive degree.
This means that the map
$\Fi{K_{S_2}+\roundup{t_2\pi^*(K_X)|_{S_2}}+L_2}$
can distinguish disjoint irreducible components in a general fiber of $\Fi{L_2}$. Thus the map
$\Fi{K_{S_2}+\roundup{(t_2+r+2)\pi^*(K_X)}|_{S_2}}$
also has this property. In order to apply (P2), we are reduced to verify
that
$$\Fi{K_{S_2}+\roundup{(t_2+r+2)\pi^*(K_X)}|_{S_2}}\bigm|_{C}$$
is an embedding for a generic irreducible element 
 $C$ of $|L_2|$. Actually, it is sufficient to verify this property for
$\Fi{K_{S_2}+\roundup{(t_2+r+2)\pi^*(K_X)|_{S_2}}}\bigm|_{C}.$

Step 2. Calculation on $C$.

We can write
$$(r+2)\pi^*(K_X)\simlin S_2+\overline{E_{\Bbb Q}},$$
where $\overline{E_{\Bbb Q}}$ is an effective ${\Bbb Q}$-divisor. So one has
$$(r+2)\pi^*(K_X)|_{S_2}\simlin S_2|_{S_2}+E_{\Bbb Q}
\simnum a_2C+E_{\Bbb Q}$$
where $E_{\Bbb Q}=\overline{E_{\Bbb Q}}|_{S_2}$ is an effective
$\Bbbq$-divisor
on $S_2$.
Considering the system 
$$\bigm|K_{S_2}+\roundup{(t_2+r+2)\pi^*(K_X)|_{S_2}
-\frac{1}{a_2}E_{\Bbb Q}}\bigm|,$$ 
we have
$$K_{S_2}+\roundup{(t_2+r+2)\pi^*(K_X)|_{S_2}}\ge
K_{S_2}+\roundup{(t_2+r+2)\pi^*(K_X)|_{S_2}
-\frac{1}{a_2}E_{\Bbb Q}}.$$
Note that
$$(t_2+r+2)\pi^*(K_X)|_{S_2}-\frac{1}{a_2}E_{\Bbb Q}-C\simnum
(t_2+r+2)\Bigl(1-\frac{r+2}{a_2(t_2+r+2)}\Bigr)\pi^*(K_X)|_{S_2},$$
which is a nef and big ${\Bbb Q}$-divisor on $S_2$ since
$\beta:=1-\frac{r+2}{a_2(t_2+r+2)}>0.$
Thus, by the vanishing theorem,
$$H^1\Bigl(S_2, K_{S_2}+ \roundup{(t_2+r+2)\pi^*(K_X)|_{S_2}-\frac{1}{a_2}E_{\Bbb
Q}}-C\Bigr)=0.$$
This would give the following exact sequence
$$H^0\Bigl(S_2,K_{S_2}+\roundup{(t_2+r+2)\pi^*(K_X)|_{S_2}-\frac{1}{a_2}E_{\Bbb
Q}}\Bigr)\lrw H^0(C, K_C+D)\lrw 0,$$
where
$D:=\roundup{(t_2+r+2)\pi^*(K_X)|_{S_2}-\frac{1}{a_2}E_{\Bbb Q}}|_C.$
Now the main task is to show that $\deg(D)\ge 3$, which implies that
$K_C+D$ is very ample since
$$\deg(K_C+D)\ge 2g(C)+1.$$
In fact, we note that $r\pi^*(K_X)|_{S_2}$ is a nef and big Cartier divisor
on $S_2$,
so \newline $r\pi^*(K_X)|_{S_2}\cdot C$ is a positive integer.
And we have
$$\align
\deg(D)&\ge \Bigl((t_2+r+2)\pi^*(K_X)|_{S_2}-\frac{1}{a_2}E_{\Bbb Q}\Bigr)\cdot C=(t_2+r+2)\beta\pi^*(K_X)\bigm|_{S_2}\cdot C\\
&=\Bigl(t_2+r+2-\frac{r+2}{a_2}\Bigr)\pi^*(K_X)\bigm|_{S_2}\cdot C\\
&\ge r\pi^*(K_X)|_{S_2}\cdot C+(3-\frac{r+2}{a_2})\pi^*(K_X)\bigm|_{S_2}\cdot C.
\endalign$$
It is easy to see $3-\frac{r+2}{a_2}>0$. So $\deg(D)\ge 3$ follows whenever 
$r\pi^*(K_X)|_{S_2}\cdot C\ge 2$, which will be proved in the next step.

Step 3. Estimating $r\pi^*(K_X)|_{S_2}\cdot C$ by studying $\fei{3r+5}.$

We claim that $r\pi^*(K_X)|_{S_2}\cdot C\ge 2$. This can be derived from our studying $\fei{3r+5}$. We have to use a lot of notations to perform the calculation. 

We know that 
$$K_{X'}+\roundup{(2r+2)\pi^*(K_X)}+S_2\le (3r+5)K_{X'}.$$
The vanishing theorem gives 
$$\align
&|K_{X'}+\roundup{(2r+2)\pi^*(K_X)}+S_2||_{S_2}=
|K_{S_2}+\roundup{(2r+2)\pi^*(K_X)}|_{S_2}|\\ 
\supset&|K_{S_2}+\roundup{(2r+2)\pi^*(K_X)|_{S_2}}\supset
|K_{S_2}+\roundup{(2r+2)\pi^*(K_X)|_{S_2}-\frac{1}{a_2}E_{\Bbb Q}}|. \tag{4.1}
\endalign$$
Because
$$(2r+2)\pi^*(K_X)|_{S_2}-\frac{1}{a_2}E_{\Bbb Q}-C\simnum (2r+2-\frac{r+2}{a_2})\pi^*(K_X)|_{S_2}$$
is a nef and big ${\Bbb Q}$-divisor, the vanishing theorem gives
$$|K_{S_2}+\roundup{(2r+2)\pi^*(K_X)|_{S_2}-\frac{1}{a_2}E_{\Bbb Q}}||_C=
|K_C+D_{3r+5}|, \tag{4.2}$$
where $D_{3r+5}:=\roundup{(2r+2)\pi^*(K_X)|_{S_2}-\frac{1}{a_2}E_{\Bbb Q}}|_C$
is a divisor on $C$ with
$$\align
\deg(D_{3r+5})&\ge (2r+2-\frac{r+2}{a_2})\pi^*(K_X)|_{S_2}\cdot C\\
&=2r\pi^*(K_X)|_{S_2}\cdot C+(2-\frac{r+2}{a_2})\pi^*(K_X)|_{S_2}\cdot C>2.
\endalign$$

Now let $M_{3r+5}$ be the movable part of $|(3r+5)K_{X'}|$ and $M_{3r+5}'$ be the movable part of $|K_{X'}+\roundup{(2r+2)\pi^*(K_X)}+S_2|$. Then it's clear
that 
$$(3r+5)\pi^*(K_X)\ge M_{3r+5}\ge M_{3r+5}'.$$
Let $L_{3r+5}$ be the movable part of $|K_{S_2}+\roundup{(2r+2)\pi^*(K_X)}|_{S_2}|$ and $L_{3r+5}'$ be the movable part of 
$$|K_{S_2}+\roundup{(2r+2)\pi^*(K_X)|_{S_2}-\frac{1}{a_2}E_{\Bbb Q}}|.$$
Then $L_{3r+5}\ge L_{3r+5}'$.  {}From (4.1) and Lemma 2.7, we have 
$M_{3r+5}|_{S_2}\ge L_{3r+5}$. {}From (4.2) and Lemma 2.7, we have
$$\align
&h^0(C, L_{3r+5}'|_C)=h^0(C, K_C+D_{3r+5})\\
=&g(C)-1+\deg(D_{3r+5})\ge g(C)+2.
\endalign$$
Using R-R and the Clifford's theorem and noting that $g(C)\ge 2$, one can easily see that 
$L_{3r+5}'\cdot C\ge 2g(C)+1\ge 5$.
So we get 
$$\align
&(3r+5)\pi^*(K_X)|_{S_2}\cdot C\ge M_{3r+5}|_{S_2}\cdot C\\
\ge& L_{3r+5}\cdot C\ge L_{3r+5}'\cdot C\ge 5
\endalign$$
i.e. $r\pi^*(K_X)|_{S_2}\cdot C\ge \frac{5r}{3r+5}>1$. Noting that 
$r\pi^*(K_X)|_{S_2}\cdot C$ is an integer, we see $r\pi^*(K_X)|_{S_2}\cdot C\ge 2$. 
The proof is completed.
\qed
\enddemo

{}From Theorem 4.1 and Theorem 4.2, we instantly have the following 
\proclaim{\smc Corollary 4.3} Let $X$ be a projective minimal 3-fold of general type
with the canonical index $r\ge 6$. Then $\fei{m_0(r)}$ is stably birational.
\endproclaim
\demo{Proof} The mian point is $\dim\fei{r+2}(X)\ge 2$ for $r\ge 6$ according to Lemma 2.4.
\qed\enddemo
%%%%%%%%%%%%%

For $4\le r\le 5$, we have to treat the case with $\dim\fei{r+2}(X)=1$. We shall use a similar method as above by studying the system $|(r+3)K_{X'}|$ because 
$\dim\fei{r+3}(X)\ge 2$. 
First we take a birational modification  $\pi: X'\lrw X$,
according to Hironaka, such that

(i) $X'$ is smooth;

(ii) the system $|(r+3)K_{X'}|$ defines a morphism;

(iii) the fractional part of $\pi^*(K_X)$ has supports with only normal
crossings.

Denote by $M_3$ the movable part of $|(r+3)K_{X'}|$. 
Take a general member $S_3\in |M_3|$. Then $S_3$ is a smooth irreducible projective surface of general type.
Set $L_3:=M_3|_{S_3}$. Then $L_3$ is a nef divisor on the surface
$S_3$. Taking $s=3$ and using (2.1) and (2.2), we have 
$$h^0(X', S_3)=P_{r+3}(X)\ge \frac{1}{8}(r^2+10r+24).$$
Thus $h^0(S_3, L_3)\ge h^0(X', S_3)-1\ge\frac{1}{8}(r^2+10r+16)$. 
%%%%%%%%%%
\proclaim{\smc Theorem 4.4}
 Let $X$ be a projective minimal 3-fold of general type
with the canonical index $r=4,\ 5$. Suppose $\dim\fei{r+3}(X)=3$. Then

(i) if $r=4$, $\fei{17}$ is stably birational;

(ii) if $r=5$, $\fei{18}$ is stably birational.
\endproclaim
\demo{Proof}
Given an integer $t_3>0$,
we note that
$$|K_{X'}+\roundup{(t_3+r+3)\pi^*(K_X)}+S_3|\subset |(t_3+2r+7)K_{X'}|.$$
In order to prove the birationality of $\fei{t_3+2r+7}$, it is sufficient to
prove the same thing for
$\Fi{K_{X'}+\roundup{(t_3+r+3)\pi^*(K_X)}+S_3}.$

Step 1. Reduce to the problem on a surface.

Since
$K_{X'}+\roundup{(t_3+r+3)\pi^*(K_X)}$ is effective by Lemma 2.5, it is enough to verify the same
thing for its restriction to $S_3$ by virtue of (P1). The vanishing theorem
gives 
$$|K_{X'}+\roundup{(t_3+r+3)\pi^*(K_X)}+S_3|\bigm|_{S_3}=
\bigm|K_{S_3}+\roundup{(t_3+r+3)\pi^*(K_X)}|_{S_3} \bigm|.$$
And from Lemma 2.2, we have
$$\bigm|K_{S_3}+\roundup{(t_3+r+3)\pi^*(K_X)|_{S_3}}\bigm|\subset
\bigm|K_{S_3}+\roundup{(t_3+r+3)\pi^*(K_X)}|_{S_3}\bigm|.$$
So, sometimes,  it is enough to show that
$\bigm|K_{S_3}+\roundup{(t_3+r+3)\pi^*(K_X)|_{S_3}}\bigm|$
gives a birational map. Under the assumption of this theorem, it is obvious
that $|L_3|$ gives a generically finite map.

Step 2. Reduce to the problem on a curve.

We suppose ${C'}$ is the general member of the movable part of $|L_3|$.
Since
$$K_{S_3}+\roundup{(t_3+r+3)\pi^*(K_X)}|_{S_3}=\bigl(K_{X'}+
\roundup{(t_3+r+3)\pi^*(K_X)}\bigr)|_{S_3}+L_3\ge L_3,$$
by (P1), we only have to verify the birationality of
$\Fi{K_{S_3}+\roundup{(t_3+r+3)\pi^*(K_X)}|_{S_3}}\bigm|_{{C'}}$
for a general member ${C'}$.
It is obvious that
$$\align
K_{S_3}+\roundup{(t_3+r+3)\pi^*(K_X)}|_{S_3}&\ge
K_{S_3}+\roundup{t_3\pi^*(K_X)}|_{S_3}+L_3\\
&\ge K_{S_3}+\roundup{t_3\pi^*(K_X)|_{S_3}}+{C'}.
\endalign $$
It's sufficient to show that $\Fi{K_{S_3}+\roundup{t_3\pi^*(K_X)|_{S_3}}+{C'}}|_{{C'}}$
is birational.

Step 3. Verifying the embedding on ${C'}$.

The vanishing theorem gives
$$\bigm|K_{S_3}+\roundup{t_3\pi^*(K_X)|_{S_3}}+{C'} \bigm|
\Bigm|_{{C'}}=
|K_{{C'}}+D_1|,$$
where
$D_1:=\roundup{t_3\pi^*(K_X)|_{S_3}}|_{{C'}}$
is a divisor on ${C'}$ with 
$$\deg(D_1)\ge t_3\pi^*(K_X)|_{S_3}\cdot {C'}.$$
It's clear that the theorem follows whenever $\deg(D_1)\ge 3$.
Although, in general, $\pi^*(K_X)|_{S_3}\cdot {C'}$ is only a rational number, we can still estimate it in this situation.

Note that, if $|{C'}|$ has already given a birational map, then 
so does $|K_{S_3}+\roundup{(t_3+r+3)\pi^*(K_X)}|_{S_3}|$ because 
$$K_{S_3}+\roundup{(t_3+r+3)\pi^*(K_X)}|_{S_3}\ge L_3.$$
So we may suppose that $|{C'}|$ gives a generically finite, non-birational map on the surface $S_3$. According to Lemma 2.6(i), we get
$${C'}^2\ge 2h^0(S_2, {C'})- 4\ge \frac{1}{4}(r^2+10r).$$
Thus we get
$$(r+3)\pi^*(K_X)|_{S_3}\cdot {C'}\ge {C'}^2\ge
\frac{1}{4}(r^2+10r).$$
One can easily obtain
$$\pi^*(K_X)|_{S_3}\cdot {C'}\cases
\ge 2,\ &\ \text{if}\ r=4\\
\ge \frac{19}{8},\ &\ \text{if}\ r=5.
\endcases$$
Thus we can see that, whenever $t_3\ge 2$ if $r=4$ or $t_3\ge 1$ if $r=5$, 
$$\deg(D_1)\ge \roundup{t_3\pi^*(K_X)|_{S_3}\cdot {C'}}\ge 3.$$
We have proved the theorem.
\qed\enddemo

%%%%%%%%%%%%%
\proclaim{\smc Theorem 4.5}
 Let $X$ be a projective minimal 3-fold of general type
with the canonical index $r=4,\ 5$. Suppose $\dim\fei{r+3}(X)=2$. Then

(i) if $r=4$, $\fei{16}$ is stably birational;

(ii) if $r=5$, $\fei{18}$ is stably birational.
\endproclaim
\demo{Proof}
In this case, we note that $|L_3|$ is a base point free pencil on the
surface
$S_3$. We can write
$$ L_3\simlin \sum_{i=1}^{a_3}C_i\simnum a_3C,$$
where
$$a_3\ge h^0(S_3, L_3)-1\ge \frac{1}{8}(r^2+10r+8)$$
and $C$ denotes a generic
irreducible element of $|L_3|$. Given an integer $t_4>0$,
we want to prove the birationality of $\fei{t_4+2r+7}$. For the same reason as in the proof of Theorem 4.4, 
it is sufficient to verify the birationality of 
$\Fi{K_{S_3}+\roundup{(t_4+r+3)\pi^*(K_X)}|_{S_3}}.$

Step 1. Reduce to the problem on a curve.

Since $K_{X'}+\roundup{(t_4+r+3)\pi^*(K_X)}$ is effective by Lemma 2.5, 
we can see that
$$K_{S_3}+\roundup{(t_4+r+3)\pi^*(K_X)}|_{S_3}\ge L_3.$$
This shows that
$\Fi{K_{S_3}+\roundup{(t_4+r+3)\pi^*(K_X)}|_{S_3}}$
can distinguish different fibers of $\Fi{L_3}$. On the other hand, we have
$$(t_4+r+3)\pi^*(K_X)\ge t_4\pi^*(K_X)+S_3.$$
So we get
$$K_{S_3}+\roundup{(t_4+r+3)\pi^*(K_X)}|_{S_3}\ge
K_{S_3}+\roundup{t_4\pi^*(K_X)}|_{S_3}+L_3.$$
The vanishing theorem gives the following exact sequence
$$H^0\Bigl(S_3, K_{S_3}+\roundup{t_4\pi^*(K_X)|_{S_3}}+L_3\Bigr)
\lrw\oplus_{i=1}^{a_3}H^0(C_i, K_{C_i}+D_i)\lrw 0$$
where the $D_i's$ are divisors on the curve $C_i$ with positive degree.
This means that the map
$\Fi{K_{S_3}+\roundup{t_4\pi^*(K_X)|_{S_3}}+L_3}$
can distinguish disjoint irreducible components in a general fiber of $\Fi{L_3}$. Thus the map
$\Fi{K_{S_3}+\roundup{(t_4+r+3)\pi^*(K_X)}|_{S_3}}$
also has this property. In order to apply (P2), we are reduced to verify
that
$$\Fi{K_{S_3}+\roundup{(t_4+r+3)\pi^*(K_X)}|_{S_3}}\bigm|_{C}$$
is an embedding for a generic irreducible element 
 $C$ of $|L_3|$. Actually, it is sufficient to verify this property for
$\Fi{K_{S_3}+\roundup{(t_4+r+3)\pi^*(K_X)|_{S_3}}}\bigm|_{C}.$

Step 2. Calculation on $C$.

We can write
$$(r+3)\pi^*(K_X)\simlin S_3+\overline{E_{\Bbb Q}},$$
where $\overline{E_{\Bbb Q}}$ is an effective ${\Bbb Q}$-divisor. So one has
$$(r+3)\pi^*(K_X)|_{S_3}\simlin S_3|_{S_3}+E_{\Bbb Q}
\simnum a_3C+E_{\Bbb Q}$$
where $E_{\Bbb Q}=\overline{E_{\Bbb Q}}|_{S_3}$ is an effective
$\Bbbq$-divisor on $S_3$.
Considering the system 
$$\bigm|K_{S_3}+\roundup{(t_4+r+3)\pi^*(K_X)|_{S_3}
-\frac{1}{a_3}E_{\Bbb Q}}\bigm|,$$ 
we have
$$K_{S_3}+\roundup{(t_4+r+3)\pi^*(K_X)|_{S_3}}\ge
K_{S_3}+\roundup{(t_4+r+3)\pi^*(K_X)|_{S_3}
-\frac{1}{a_3}E_{\Bbb Q}}.$$
Note that
$$(t_4+r+3)\pi^*(K_X)|_{S_3}-\frac{1}{a_3}E_{\Bbb Q}-C\simnum
(t_4+r+3)\Bigl(1-\frac{r+3}{a_3(t_4+r+3)}\Bigr)\pi^*(K_X)|_{S_3},$$
which is a nef and big ${\Bbb Q}$-divisor on $S_3$ since
$\gamma:=1-\frac{r+3}{a_3(t_4+r+3)}>0.$
Thus, by the vanishing theorem,
$$H^1\Bigl(S_3, K_{S_3}+ \roundup{(t_4+r+3)\pi^*(K_X)|_{S_3}-\frac{1}{a_3}E_{\Bbb Q}}-C\Bigr)=0.$$
This would give the following exact sequence
$$H^0\Bigl(S_3,K_{S_3}+\roundup{(t_4+r+3)\pi^*(K_X)|_{S_3}-\frac{1}{a_3}E_{\Bbb
Q}}\Bigr)\lrw H^0(C, K_C+D)\lrw 0,$$
where
$D:=\roundup{(t_4+r+3)\pi^*(K_X)|_{S_3}-\frac{1}{a_3}E_{\Bbb Q}}|_C.$
Now the main task is to show that $\deg(D)\ge 3$, which implies that
$K_C+D$ is very ample since
$$\deg(K_C+D)\ge 2g(C)+1.$$
In fact, we note that $r\pi^*(K_X)|_{S_3}$ is a nef and big Cartier divisor
on $S_3$,
so \newline $r\pi^*(K_X)|_{S_3}\cdot C$ is a positive integer.
And we have
$$\align
\deg(D)&\ge \Bigl((t_4+r+3)\pi^*(K_X)|_{S_3}-\frac{1}{a_3}E_{\Bbb Q}\Bigr)\cdot C=(t_4+r+3)\gamma\pi^*(K_X)\bigm|_{S_3}\cdot C\\
&=\Bigl(t_4+r+3-\frac{r+3}{a_3}\Bigr)\pi^*(K_X)\bigm|_{S_3}\cdot C\\
&\ge r\pi^*(K_X)|_{S_3}\cdot C+(4-\frac{r+3}{a_3})\pi^*(K_X)\bigm|_{S_3}\cdot C.
\endalign$$
It is easy to see $4-\frac{r+3}{a_3}>0$. So $\deg(D)\ge 3$ follows whenever 
$r\pi^*(K_X)|_{S_3}\cdot C\ge 2$, which will be proved in the next step.

Step 3. Estimating $r\pi^*(K_X)|_{S_3}\cdot C$ by studying $\fei{3r+5}.$

We claim that $r\pi^*(K_X)|_{S_3}\cdot C\ge 2$. This can be derived from our studying $\fei{3r+5}$. We have to use a lot of notations to perform the calculation. 

We know that 
$$K_{X'}+\roundup{(2r+1)\pi^*(K_X)}+S_3\le (3r+5)K_{X'}.$$
The vanishing theorem gives 
$$\align
&|K_{X'}+\roundup{(2r+1)\pi^*(K_X)}+S_3||_{S_3}=
|K_{S_3}+\roundup{(2r+1)\pi^*(K_X)}|_{S_3}|\\ 
\supset&|K_{S_3}+\roundup{(2r+1)\pi^*(K_X)|_{S_3}}\supset
|K_{S_3}+\roundup{(2r+1)\pi^*(K_X)|_{S_3}-\frac{1}{a_3}E_{\Bbb Q}}|. \tag{4.3}
\endalign$$
Because
$$(2r+1)\pi^*(K_X)|_{S_3}-\frac{1}{a_3}E_{\Bbb Q}-C\simnum (2r+1-\frac{r+3}{a_3})\pi^*(K_X)|_{S_3}$$
is a nef and big ${\Bbb Q}$-divisor, the vanishing theorem gives
$$|K_{S_3}+\roundup{(2r+1)\pi^*(K_X)|_{S_3}-\frac{1}{a_3}E_{\Bbb Q}}||_C=
|K_C+D_{3r+5}|, \tag{4.4}$$
where $D_{3r+5}:=\roundup{(2r+1)\pi^*(K_X)|_{S_3}-\frac{1}{a_3}E_{\Bbb Q}}|_C$
is a divisor on $C$ with
$$\align
\deg(D_{3r+5})&\ge (2r+1-\frac{r+3}{a_3})\pi^*(K_X)|_{S_3}\cdot C\\
&=2r\pi^*(K_X)|_{S_3}\cdot C+(1-\frac{r+3}{a_3})\pi^*(K_X)|_{S_3}\cdot C>2.
\endalign$$

Now let $M_{3r+5}$ be the movable part of $|(3r+5)K_{X'}|$ and $M_{3r+5}'$ be the movable part of $|K_{X'}+\roundup{(2r+1)\pi^*(K_X)}+S_3|$. Then it's clear
that 
$$(3r+5)\pi^*(K_X)\ge M_{3r+5}\ge M_{3r+5}'.$$
Let $L_{3r+5}$ be the movable part of $|K_{S_3}+\roundup{(2r+1)\pi^*(K_X)}|_{S_3}|$ and $L_{3r+5}'$ be the movable part of 
$$|K_{S_3}+\roundup{(2r+1)\pi^*(K_X)|_{S_3}-\frac{1}{a_3}E_{\Bbb Q}}|.$$
Then $L_{3r+5}\ge L_{3r+5}'$.  {}From (4.3) and Lemma 2.7, we have 
$M_{3r+5}|_{S_3}\ge L_{3r+5}$. {}From (4.4) and Lemma 2.7, we have
$$\align
&h^0(C, L_{3r+5}'|_C)=h^0(C, K_C+D_{3r+5})\\
=&g(C)-1+\deg(D_{3r+5})\ge g(C)+2.
\endalign$$
Using R-R and the Clifford's theorem and noting that $g(C)\ge 2$, one can easily see that 
$L_{3r+5}'\cdot C\ge 2g(C)+1\ge 5$.
So we get 
$$\align
&(3r+5)\pi^*(K_X)|_{S_3}\cdot C\ge M_{3r+5}|_{S_3}\cdot C\\
\ge& L_{3r+5}\cdot C\ge L_{3r+5}'\cdot C\ge 5
\endalign$$
i.e. $r\pi^*(K_X)|_{S_3}\cdot C\ge \frac{5r}{3r+5}>1$. Noting that 
$r\pi^*(K_X)|_{S_3}\cdot C$ is an integer, we see $r\pi^*(K_X)|_{S_3}\cdot C\ge 2$. 
The proof is completed.
\qed\enddemo

{}From Theorem 3.1, we can take $m_0(2)=11$ and $m_0(3)=15$. Theorems 4.4 and 4.5 imply $m_0(4)=17$ and $m_0(5)=18$. Therefore the main theorem follows.

\head {\rm 5. Threefolds with positive geometric genus}\endhead
%%%%%%%%%%%%%%%%%

Throughout this section, we still suppose $X$ is a projective minimal 3-fold
of general type. Our aim is to study 3-folds with big geometric genus using
the method of the Main Theorem. Koll\'ar (Corollary 4.8 of \cite{11}) proved
that $\fei{16}$
is birational if $p_g(X)\ge 2$. Reviewing the parallel results on surfaces
and
Gorenstein 3-folds, one should expect a better bound for the birationality of
$\fei{m}$.
\bigskip
To begin the argument, we first take a birational modification $\pi: X'\lrw
X$ according to Hironaka such that

(i) $X'$ is smooth;

(ii) $|K_{X'}|$ gives a morphism;

(iii) the fractional part of $\pi^*(K_X)$ has supports with only normal
crossings.

Set $g:=\fei{1}\circ \pi$ and take the Stein-factorization of
$$g: X'\overset{f}\to\lrw W\lrw W'\subset{\Bbb P}^N$$
where $W'$ is the image
of $X'$ through $g$ and $f$ is a fibration. Let $M$ be the movable part of
$|K_{X'}|$. We can write
$$K_{X'}\simlin M+E'\ \ \text{and}\ \ \pi^*(K_X)=_{\Bbb Q} M+\overline{E_{\Bbb Q}}$$
where $E'$ is an effective divisor and $\overline{E_{\Bbb Q}}$ is an
effective
${\Bbb Q}$-divisor.
%%%%%%%%%%%%

\proclaim{\smc Proposition 5.1}  Let $X$ be a minimal projective threefold of
general
type. If $\dim\fei{1}(X)\ge 2$, then

(i) $\fei{4}$ is generically finite;

(ii) $\fei{3}$ is generically finite provided $p_g(X)\ge 4$.
\endproclaim
\demo{Proof}
Because $p_g(X)>0$, it is sufficient to prove for the case when
$\dim\fei{1}(X)=2$. Let $S\in |M|$ be the general member. Then $S$ is a
smooth projective
surface of general type. We have $S|_S\simnum aC$, where $C$ is a smooth
curve
and $a\ge p_g(X)-2$.
Considering the system $|K_{X'}+\roundup{\pi^*(K_X)}+2S|$, we have
$$|K_{X'}+\roundup{\pi^*(K_X)}+2S|\bigm|_S=\bigm|K_S+\roundup{\pi^*(K_X)
}|_S+
S|_S\bigm|.$$
Besides, we have
$$\bigm|K_S+\roundup{\pi^*(K_X)|_S}+S|_S\bigm|\Bigm|_C=|K_C+D|,$$
where $D:=\roundup{\pi^*(K_X)|_S}\bigm|_C$ is a divisor on $C$ of positive
degree.
Thus $|K_C+D|$ gives a finite map and so does $\fei{4}$. This derives (i).
If $p_g(X)\ge 4$, then $a\ge 2$. By the vanishing theorem, we get
$$\align
|K_{X'}+\roundup{\pi^*(K_X)}+S|\bigm|_S&=|K_S+\roundup{\pi^*(K_X)}|_S|\\
&\supset |K_S+\roundup{\pi^*(K_X)|_S}|.
\endalign$$
We can write
$\pi^*(K_X)|_S=S|_S+E_{\Bbbq},$
where $E_{\Bbbq}$ is an effective $\Bbbq$-divisor. It is obvious that
$$\pi^*(K_X)|_S-\frac{1}{a}E_{\Bbbq}-C\simnum
\Bigl(1-\frac{1}{a}\Bigr)\pi^*(K_X)|_S.$$
So
$$H^1\Bigl(S, K_S+\roundup{\pi^*(K_X)|_S-\frac{1}{a}E_{\Bbbq}}-C\Bigr)=0,$$
which gives
$$|K_S+\roundup{\pi^*(K_X)|_S-\frac{1}{a}E_{\Bbbq}}|\bigm|_C=|K_C+D|,$$
where $D:=\roundup{\pi^*(K_X)-\frac{1}{a}E_{\Bbbq}}|_C$ is a divisor on
$C$ of
positive degree. Thus $|K_C+D|$ gives a generically finite map. By (P2),
$\fei{3}$
is also generically finite.
\qed
\enddemo
\proclaim{\smc Proposition 5.2}  Let $X$ be a minimal projective threefold of
general type.
If $\dim\fei{1}(X)\ge 2$, then

(i) $\fei{8}$ is birational;

(ii) $\fei{6}$ is birational provided $p_g(X)\ge 4$.
\endproclaim
\demo{Proof}
If $\dim\fei{1}(X)=3$, then it is very easy to prove the birationality of
$\fei{6}$ by standard argument.
We mainly discuss the case when $\dim\fei{1}(X)=2$.
To prove (i), let $M_4$ be the movable part of $|4K_{X'}|$. We can modify
$\pi$,
if necessary, such that $|M_4|$ is also base point free. We have
$$|K_{X'}+\roundup{\pi^*(K_X)}+M_4+2S|\bigm|_S=\bigm|K_S+\roundup{\pi^*(
K_X)}|_S+
L_4+S|_S\bigm|,$$
where $L_4:=M_4|_S$ is nef and $\Fi{L_4}$ is generically finite. It is not
difficult
to see that the right system gives a birational map using the method which
has been
applied frequently in this paper. So $\fei{8}$ is birational according to (P1).

If $p_g(X)\ge 4$,
denote by $M_3$ the moving part of $|3K_{X'}|$. For the same reason, we can
suppose
$|M_3|$ is also base point free. We have
$$\align
&|K_{X'}+\roundup{4\pi^*(K_X)}+S|\bigm|_S=|K_S+\roundup{4\pi^*(K_X)}|_S|
\\
&\supset |K_S+\roundup{4\pi^*(K_X)|_S}|\supset
|K_S+\roundup{\pi^*(K_X)|_S}
+L_3|.
\endalign$$
where $L_3:=M_3|_S$, which is a nef and big divisor on $S$. And $|L_3|$
gives
a generically finite map. Because
$$\pi^*(K_X)|_S-C-\frac{1}{a}E_{\Bbbq}$$
is nef and big, the vanishing theorem will imply that
$$|K_S+\roundup{\pi^*(K_X)|_S-\frac{1}{a}E_{\Bbbq}}+L_3|\bigm|_C=|K_C+L_3|_C+
D|,$$
where $D=\roundup{\pi^*(K_X)|_S-\frac{1}{a}E_{\Bbbq}}|_C$ is a divisor of
positive
degree. Thus $\deg(L_3|_C+D)\ge 3$ and $K_C+L_3|_C+D$ is very ample. This
shows
that  $\fei{6}$ is birational.
\qed
\enddemo
%%%%%%%%%%%

\proclaim{\smc Proposition 5.3} Let $X$ be a minimal projective threefold of
general type. If $\dim\fei{1}(X)=1$, then

(i) $\fei{9}$ is birational;

(ii) $\fei{6}$ is birational provided $p_g(X)\ge 12$.
\endproclaim
\demo{Proof} (i).
In this case, $W$ is a nonsingular curve. We set $b:=g(W)$, the genus of
$W$.

If $b>0$, then $\fei{1}$ is actually a morphism. In this case, there is no
need
to make the modification $\pi$, i.e. $X'=X$. Though $K_X$ is not Cartier, it
is
a Weil divisor. We can still define the system $|K_X|$ in a natural way.
We have $M\simlin \sum_{i=1}^a S_i$, where the $S_i$ are fibers of $f$.
Noting that
the singularities on $X$ are all isolated, a general $S_i$ is a smooth
projective
surface of general type. Using Kawamata's vanishing theorem (\cite{10})
for $\Bbbq$-Cartier Weil divisors, we have $H^1(X,\omega_X^{[k]})=0$ whenever
$k>1$. Thus \newline
$\fei{6}\bigm|_{S_i}=\Fi{5K_{S_i}}$
is birational. According to (P2), $\fei{6}$ is birational.

If $b=0$, $W=\Bbbp1$ and we have the fibration $f:X'\lrw \Bbbp1$. Let $S$ be
a general
fiber of $f$. Then $S$ is a smooth projective surface of general type. We
divide $S$ into two categories:

(I) $(K_{S_0}^2, p_g(S))\ne (1,2)$ and $(2,3)$;

(II) the rest.

\noindent where $S_0$ denotes the minimal model of $S$. 

Suppose $S$ is of type (I). There is a common property for these surfaces
that the 3-canonical maps are birational. To deal with this situation,
We can use Koll\'ar's technique (the proof of Corollary 4.8 in \cite{11}).  Because $p_g(X)>0$, we have $p_g(S)>0$.  Let $\sigma: S\lrw S_0$ be the contraction onto the minimal model. 
According to Theorem 3.1 in \cite{5}, we see that $|2K_{S_0}|$ is base point free.  So the movable part of $|2K_S|$ is $\sigma^*(2K_{S_0})$.
We have
$H^0(\omega_{X'}^{7})=H^0({\Bbb P}^1, f_*\omega_{X'}^{7})$
 and an
 injection $\Cal{O}(1)\hookrightarrow f_*\omega_{X'}$, and hence an
injection
$\Cal{O}(5)\hookrightarrow f_*\omega_{X'}^{5}.$
This gives an injection
$$\Cal{O}(5)\otimes f_*\omega_{X'}^2\hookrightarrow f_*\omega_{X'}^{7},$$
where
$$\Cal{O}(5)\otimes f_*\omega_{X'}^2=
\Cal{O}(1)\otimes f_*\omega_{X'/{\Bbb P}^1}^2.$$
It is well-known that $f_*\omega_{X'/{\Bbb P}^1}^2$ is a sum of line bundles
of non-negative degree on ${\Bbb P}^1$. The local sections of
$f_*\omega_{X'}^2$ give the bicanonical map for $S$, and all these extend to
global sections of $\Cal{O}(5)\otimes f_*\omega_{X'}^2$. Moreover the sections of $\Cal{O}(1)\otimes f_*\omega_{X'/{\Bbb P}^1}^2$ 
separate different fibers. 
Suppose $M_7$ is the movable part of $|7K_{X'}|$. 
Because $\fei{7}=\Fi{M_7}$, 
we can see from the above argument that 
$M_7|_S\ge \sigma^*(2K_{S_0}).$
Now considering the system 
$|K_{X'}+\roundup{7\pi^*(K_X)}+S|,$
we have 
$$\align
&|K_{X'}+\roundup{7\pi^*(K_X)}+S|\bigm|_S=\bigm|K_S+\roundup{7\pi^*(K_X)}|_S
\bigm|\\
&\supset |K_S+M_7|_S|\supset |K_S+\sigma^*(2K_{S_0})|.
\endalign$$
Because $|K_S+\sigma^*(2K_{S_0})|$ gives a birational map, 
we see that $\fei{9}$ is birational.

Suppose $S$ is of type (II) and $(K_{S_0}^2,p_g(S))=(2,3)$. We want to show that $\fei{8}$ is birational. We know that the movable part of $|K_{X'}|$ is linearly equivalent to a disjoint union of irreducible copies of $S$.
We  have
$$|K_{X'}+\roundup{\pi^*(K_X)}+S|\bigm|_S=|K_S+
\roundup{\pi^*(K_X)}|_S|\supset |K_S|.$$
Because the movable part of $|K_S|$ gives a finite map onto ${\Bbb P}^2$, we can see that $\fei{3}$ is generically finite.
Denote by $M_3$ the movable part of $|3K_{X'}|$ and by $M_5$ the movable part of
$|5K_{X'}|$. In order to prove the birationality of $\fei{8}$, we should
study $\bigm|M_5|_S\bigm|.$ 
We have
$K_{X'}+\roundup{3\pi^*(K_X)}+S\le 5K_{X'}.$
The vanishing theorem gives 
$$|K_{X'}+\roundup{3\pi^*(K_X)}+S|\bigm|_S=|K_S+\roundup{3\pi^*(K_X)}|_S|.$$
We suppose that $K_0$ is the movable part of $|K_S|$.
Denote by $M_{+}$ the movable part of
$|K_{X'}+\roundup{3\pi^*(K_X)}+S|$. Then $M_{+}\le M_5$. 
By Lemma 2.7, we can see that ${M_{+}}|_S$ contains the movable part of 
$|K_S+\roundup{3\pi^*(K_X)}|_S|$. Let $L$ be the movable part of $|{M_3}|_S|$. Then
$K_S+\roundup{3\pi^*(K_X)}|_S\ge K_0+L$
and $K_0+L$ is movable. So $M_{+}|_S\ge K_0+L.$
On the other hand,
because $M_5\ge M_{+}$, so $M_5|_S\ge K_0+L$. Now it is time to study the
$\fei{8}$.
For an obvious reason, we can suppose that $|M_5|$ is base point free.
This
assumption means that $M_5$ is nef and big. The vanishing theorem gives the
exact sequence
$$H^0\Bigl(X', K_{X'}+\roundup{\pi^*(K_X)}+M_5+S\Bigr)\lrw
H^0\bigl(S, K_S+\roundup{\pi^*(K_X)}|_S+M_5|_S\bigr)\lrw 0.$$
We note that
$$K_S+\roundup{\pi^*(K_X)}|_S+M_5|_S\ge
K_S+\roundup{\pi^*(K_X)|_S}+L+K_0.$$
Since $\pi^*(K_X)|_S$ is nef and big, $\roundup{\pi^*(K_X)|_S}$ is
effective, $\dim\Fi{K_0}(S)=2$ and
$\dim\Fi{L}(S)=2$, using our method again, it is easy to see that
$|K_S+\roundup{\pi^*(K_X)|_S}+L+K_0|$ gives a birational map. Which shows
that
$$|K_{X'}+\roundup{\pi^*(K_X)}+M_5+S|$$
gives a birational map and so does $\fei{8}$. 

Suppose $S$ is of type (II) and $(K_{S_0}^2, p_g(S))=(1,2)$. We want to show that $\fei{9}$ is birational. This is the most frustrating case because $\Fi{4K_S}$ is not birational. We recall that $|K_{S_0}|$ has no fixed component, that it has exactly one base point and that a general member of this system is a smooth irreducible curve of genus $2$.  Thus the movable part $C$ of $|K_S|$ is also a smooth curve of genus two. Furthermore $C\le \sigma^*(K_{S_0})$. Because 
$$|K_{X'}+\roundup{\pi^*(K_X)}+S|\bigm|_S
=|K_S+\roundup{\pi^*(K_X)}|_S|\supset |K_S|\supset |C|,$$
we see that $\dim\fei{3}(X)\ge 2$. We still denote by $M_3$ the movable part of $|3K_{X'}|$ and by $M_3'$ the movable part of $|K_{X'}+\roundup{\pi^*(K_X)}+S|$.
According to Lemma 2.7, $M_3|_S\ge M_3'|_S\ge C$. Now we consider the system
$$|K_{X'}+\roundup{4\pi^*(K_X)}+M_3+S|.$$
Actually we  can take further modification to $\pi$ such that $|M_3|$ is also base point free. This means we can suppose $M_3$ is nef. By the Kawamata-Viehweg vanishing theorem, we have
$$\align
&|K_{X'}+\roundup{4\pi^*(K_X)}+M_3+S|\bigm|_S
=\bigm|K_S+\roundup{4\pi^*(K_X)}|_S+M_3|_S\bigm|\\
\supset&\bigm|K_S+\roundup{4\pi^*(K_X)|_S}+C\bigm|.
\endalign$$
We can use the vanishing theorem once more so that we get
$$\bigm|K_S+\roundup{4\pi^*(K_X)|_S}+C\bigm|\bigm|_C=|K_C+D|,$$
where $D:=\roundup{4\pi^*(K_X)|_S}\bigm|_C.$
So if we can prove $4\pi^*(K_X)|_S\cdot C>2$, then $\deg(D)\ge 3$ which induces the birationality of $\fei{9}$, because
$$|K_{X'}+\roundup{4\pi^*(K_X)}+M_3+S|\subset |9K_{X'}|.$$
Thus we only have to prove the following claim.
\medskip

\noindent {\bf Claim.}\ \  $\xi:=\pi^*(K_X)|_S\cdot C\ge \frac{3}{5}$.
\smallskip

The idea is to find an initial estimation to $\xi$ by first proving that $\fei{10}$ 
is birational. Then we can optimize this estimation by an infinite programme. We will find that the limit estimation is $\frac{3}{5}$.  Actually our second estimation is enough for us to show the birationality of $\fei{9}$. We present a better estimation here hoping that it might be useful to prove the birationality of $\fei{8}$ in future.  

From now on, we prove the claim. Let $M_5$, $M_7$ and $M_{10}$ be the movable part of $|5K_{X'}|$, $|7K_{X'}|$ and $|10K_{X'}|$ respectively. For the same reason, we can suppose they are all nef. 
We can see that
$$\align
&|K_{X'}+\roundup{3\pi^*(K_X)}+S|\bigm|_S
=\bigm|K_S+\roundup{3\pi^*(K_X)}|_S\bigm|\\
\supset&\bigm|K_S+M_3|_S\bigm|\supset|K_S+C|.
\endalign$$
Because $q(S)=0$, $|K_S+C|$ gives a generically finite map. So $\fei{5}$ is generically finite. Suppose $L_5$ is the movable part of $\bigm|M_5|_S\bigm|.$ Then $\dim\Fi{L_5}(S)=2$. Therefore we can see that $L_5\cdot C\ge 2$ for a general element $C$.
We also have 
$$\align
&|K_{X'}+\roundup{5\pi^*(K_X)}+S|\bigm|_S\\
=&\bigm|K_S+\roundup{5\pi^*(K_X)}|_S\bigm|
\supset |C+L_5|. \tag{5.1}
\endalign$$
Noting that $C+L_5$ is movable and by Lemma 2.7, we see that 
$M_7'|_S\ge C+L_5$ where $M_7'$ is the movable part of 
$|K_{X'}+\roundup{5\pi^*(K_X)}+S|.$
Because $M_7\ge M_7'$, $M_7|_S\ge C+L_5$. 
We have
$$\align
&|K_{X'}+\roundup{\pi^*(K_X)}+M_7+S|\bigm|_S
=\bigm|K_S+\roundup{\pi^*(K_X)}|_S+M_7|_S\bigm|\\
\supset &\bigm|K_S+\roundup{\pi^*(K_X)|_S}+L_5+C\bigm|. \tag{5.2}
\endalign$$
Now it is obvious that 
$$
\bigm|K_S+\roundup{\pi^*(K_X)|_S}+L_5+C\bigm|\bigm|_C
=|K_C+G|,$$
where $G:=(\roundup{\pi^*(K_X)|_S}+L_5)|_C$ is a divisor of degree $\ge 3$ and so $h^0(C, K_C+G)\ge g(C)+2=4$.
Suppose $L_{10}$ is the movable part of $\bigm|M_{10}|_S\bigm|$ and 
$M_{10}'$ is the movable part of 
$|K_{X'}+\roundup{\pi^*(K_X)}+M_7+S|.$
Let $L_{10}'$ be the movable part of 
$$|K_S+\roundup{\pi^*(K_X)|_S}+L_5+C|.$$
By Lemma 2.7, we have $M_{10}'|_S\ge L_{10}'$ and 
$$h^0(C, L_{10}'|_C)=h^0(K_C+G)\ge 4.$$
Noting that $M_{10}\ge M_{10}'$, we have $L_{10}\ge L_{10}'$ and 
$h^0(C, L_{10}|_C)\ge 4.$
Because $C$ is a curve of genus $2$, by R-R, we see that $L_{10}\cdot C\ge 5$.
This means that 
$$10\pi^*(K_X)|_S\cdot C\ge M_{10}|_S\cdot C\ge  L_{10}\cdot C\ge 5.$$
So we get $\xi\ge \frac{1}{2}$.

Suppose $M_{12}$ is the movable part of $|12K_{X'}|$.
Similar to (5.1), we have
$$\align
&|K_{X'}+\roundup{10\pi^*(K_X)}+S||_S=|K_S+\roundup{10\pi^*(K_X)}|_S|\\
\supset& |K_S+L_{10}|\supset |C+L_{10}|.
\endalign$$
Using Lemma 2.7, we can easily see that $M_{12}|_S\ge C+L_{10}$. Replacing 
$M_7$ by $M_{12}$ in (5.2), we also have 
$$\align
&|K_{X'}+\roundup{\pi^*(K_X)}+M_{12}+S||_S=|K_S+\roundup{\pi^*(K_X)}|_S
+M_{12}|_S|\\
\supset& |K_S+\roundup{\pi^*(K_X)|_S}+L_{10}+C|.
\endalign$$
Using the vanishing theorem once more, we have
$$|K_S+\roundup{\pi^*(K_X)|_S}+L_{10}+C||_C=|K_C+\roundup{\pi^*(K_X)|_S}|_C+
L_{10}|_C|.$$
Let $L_{15}'$ be the movable part of $|K_S+\roundup{\pi^*(K_X)|_S}+L_{10}+C|$.
{}From Lemma 2.7, it is easy to see $M_{15}|_S\ge L_{15}'$ and
$$\align
&h^0(C, L_{15}'|_C)=h^0(K_C+\roundup{\pi^*(K_X)|_S}|_C+L_{10}|_C)\\
=&\deg(\roundup{\pi^*(K_X)|_S}|_C)+\deg(L_{10}|_C)+g(C)-1\ge 7
\endalign$$
where $M_{15}$ is the movable part of $|15K_{X'}|$. 
By R-R and the Clifford's theorem, we see that
$$L_{15}'\cdot C\ge h^0(C, L_{15}'|_C)+g(C)-1\ge 8.$$
Thus
$$15\pi^*(K_X)|_S\cdot C\ge M_{15}|_S\cdot C\ge L_{15}'\cdot C\ge 8.$$
This means that 
$\xi\ge \frac{8}{15}>\frac{1}{2},$ which directly induces the birationality of 
$\fei{9}$.

We can infinitely repeat this programme, but omit the details.  
So we can get the following sequence
$$\align
&n_0=10,\ d_0=5\\
&n_1=n_0+5,\ d_1=d_0+3\\
&\cdots\cdots\cdots\\
&n_k=n_{k-1}+5,\ d_k=d_{k-1}+3\\
&\cdots\cdots\cdots\\
&\xi\ge\frac{d_k}{n_k},\ \text{for all}\ k.
\endalign$$
Thus 
$$\xi\ge \lim_{k\mapsto\infty}\frac{d_k}{n_k}=
\lim_{k\mapsto\infty}\frac{3k+5}{5k+10}=\frac{3}{5}.$$
The claim is proved.

(ii). If $p_g(X)\ge 12$, then we have
${\Cal O}(11)\hookrightarrow f_*\omega_{X'}$. So
$${\Cal O}(1)\otimes f_*\omega_{X'/\Bbbp1}^5={\Cal O}(11)\otimes
f_*\omega_{X'}^5
\hookrightarrow f_*\omega_{X'}^6.$$
It is easy to see that $\fei{6}$ is birational for $X$ by virtue of
Koll\'ar's technique.
\qed
\enddemo

Both Proposition 5.2 and Proposition 5.3 imply
\proclaim{\smc Corollary 5.4} Let $X$ be a smooth projective 3-fold of general
type.
Then $\fei{9}$ is birational if $p_g(X)\ge 2$. 
\endproclaim
\medskip

\noindent{\bf Acknowledgment.}
This paper was begun while I was visiting the Abdus Salam International Centre for Theoretical Physics, Trieste, Italy in 1999 and was finally fulfilled 
when I stayed as a visiting fellow at the Mathematisches Institut der Universit$\ddot{\text{a}}$t G$\ddot{\text{o}}$ttingen, Germany in 2000. I would like to thank F. Catanese for fruitful discussions during the preparation of it. Special thanks are due to hospital faculty  members of the Mathematisches Institut der Universit$\ddot{\text{a}}$t G$\ddot{\text{o}}$ttingen. 

\head {\bf References} \endhead
\roster
\item"[1]" E. Bombieri, Canonical models of surfaces of general type,
{\it Publications I.H.E.S.} {\bf 42}(1973), 171-219.
\item"[2]" F. Catanese, Canonical rings and special surfaces of
general type, {\it Proc. Symp. Pure Math.} {\bf 46}(1987), 175-194.
\item"[3]" M. Chen, On pluricanonical maps for threefolds of general
type, {\it J. Math. Soc. Japan} {\bf 50}(1998), 615-621.
\item"[4]" -----, The relative pluricanonical stability for 3-folds
of general type, to appear in {\it Proc. Amer. Math. Soc.}
\item"[5]" C. Ciliberto, The bicanonical map for surfaces of general type, {\it Proc. Symposia in Pure Math.} {\bf 62}(1997), 57-83.
\item"[6]" L. Ein, R. Lazarsfeld, Global generation of pluricanonical
and adjoint linear systems on smooth projective threefolds, {\it J. Amer. Math. Soc.} {\bf 6}(1993), 875-903.
\item"[7]" M. Hanamura, Stability of the pluricanonical maps of
threefolds, {\it Adv. Stud. Pure Math.} {\bf 10}(1987), 185-205.
\item"[8]" R. Hartshorne, {\it Algebraic Geometry}, GTM {\bf 52},
Springer-Verlag
1977.
\item"[9]" Y. Kawamata, A generalization of Kodaira-Ramanujam's
vanishing theorem, {\it Math. Ann.} {\bf 261}(1982), 43-46.
\item"[10]" Y. Kawamata, K. Matsuda, K. Matsuki, Introduction to the
minimal model problem, {\it Adv. Stud. Pure Math.} {\bf 10}(1987), 283-360.
\item"[11]" J. Koll\'ar, Higher direct images of dualizing sheaves, I,
{\it Ann. of Math.} {\bf 123}(1986), 11-42.
\item"[12]" J. Koll\'ar, S. Mori, {\it Birational geometry of algebraic
varieties}, Cambridge Univ. Press, 1998.
\item"[13]" S. Lee, Remarks on the pluricanonical and adjoint linear
series on projective threefolds, {\it Commun. Algebra} {\bf 27}(1999), 4459-4476.
\item"[14]" K. Matsuki, On pluricanonical maps for 3-folds of general
type, {\it J. Math. Soc. Japan} {\bf 38}(1986), 339-359.
\item"[15]" M. Reid, Canonical 3-folds, {\it Journ\'ees de G\'eom\'etrie
Alg\'ebrique d'Angers} (Beauville Ed.), Sijthoff and Noordhof, Alphen ann den
Rijn 1980, 273-310.
\item"[16]" I. Reider, Vector bundles of rank 2 and linear systems on
algebraic surfaces, {\it Ann. Math.} {\bf 127}(1988), 309-316.
\item"[17]" S. G. Tankeev, On n-dimensional canonically polarized
varieties and varieties of fundamental type, {\it Izv. A. N. SSSR}, S\'er. Math. {\bf 35}(1971), 31-44.
\item"[18]" E. Viehweg, Vanishing theorems, {\it J. reine angew. Math.} {\bf
335}(1982),
1-8.
\item"[19]" P.M.H. Wilson, Towards birational classification of algebraic varieties, {\it Bull. Lond. Math. Soc.} {\bf 19}(1987), 1-48. 
\item"[20]" ------, The pluricanonical map on varieties of general type, {\it Bull. Lond. Math. Soc.} {\bf 12}(1980), 103-107.
\item"[21]" G. Xiao, Linear bound for abelian automorphisms of varieties
of general type, {\it J. reine angew. Math.} {\bf 476}(1996), 201-207.
\endroster
\enddocument